\documentclass[12pt]{amsart}

\newtheorem{thm}{Theorem}[section]

\newtheorem{lem}{Lemma}[section]

\theoremstyle{definition}
\newtheorem{df}{Definition}[section]

\theoremstyle{remark}

\numberwithin{equation}{section}

\begin{document}
\title[On $\phi$-symmetric LP-Sasakian Manifolds ....]{On $\phi$-symmetric LP-Sasakian Manifolds
admitting semi-symmetric metric connection}
\author[A. A. Shaikh and S. K. Hui]{Absos Ali Shaikh and Shyamal Kumar Hui}

\begin{abstract}
The object of the present paper is to study locally $\phi$-symmetric LP-Sasakian manifolds admitting
semi-symmetric metric connection and obtain a necessary and sufficient condition for a locally $\phi$-symmetric LP-Sasakian
manifold with respect to semi-symmetric metric connection to be locally $\phi$-symmetric LP-Sasakian
manifold with respect to Levi-Civita connection.
\end{abstract}
\maketitle
\section{Introduction}
On the analogy of Sasakian manifolds, in 1989 Matsumoto \cite{MAT} introduced the notion of LP-Sasakian
manifolds. Again the same notion was studied by Mihai and Rosca \cite{MIHAI2} and they obtained many results.
LP-Sasakian manifolds are also studied by De et. al. \cite{DE0}, Shaikh et. al. (\cite{SHAIKH2},
\cite{SHAIKH1}, \cite{SHAIKH3}, \cite{SHAIKH5}), Taleshian and Asghari \cite{TAL}, Venkatesha
and Bagewadi \cite{BAGE2} and many others. The notion of locally $\phi$-symmetry on a $3$-dimensional
LP-Sasakian manifold was studied by Shaikh and De \cite{SHAIKH8}.\\
\indent In 1924 Friedmann and Schouten \cite{FRIED} introduced the notion of semi-symmetric linear
connection on a differentiable manifold. Then in 1932 Hayden \cite{HAYDEN} introduced the idea of
metric connection with torsion on a Riemannian manifold. A systematic study of the semi-symmetric metric
connection on a Riemannian manifold has been given by Yano \cite{YANO5} in 1970. Also semi-symmetric metric
connection on a Riemannian manifold has been studied by Barua and Mukhopadhyay \cite{BARUA1}, Binh \cite{BINH1},
Chaki and Chaki \cite{cha2}, Chaturvedi and Pandey \cite{CHATUR}, Shaikh and Hui \cite{HUI4}, Sharfuddin and
Hussain \cite{sh} and many others. Recently Shaikh and Jana studied the quarter-symmetric metric
connection on a $(k,\mu)$-contact metric manifold \cite{SHAIKH9}.\\
\indent The study of Riemann symmetric manifolds began with the work of Cartan \cite{ca}.
A Riemannian manifold $(M^{n},g)$ is said to be locally symmetric due to Cartan \cite{ca}
if its curvature tensor $R$ satisfies the relation $\nabla R =0$, where $\nabla$ denotes the
operator of covariant differentiation with respect to the metric tensor $g$.
As a weaker version of local symmetry, the notion of locally $\phi$-symmetric Sasakian manifold was
introduced by Takahashi \cite{TAKA}. In the sense of Takahashi, Shaikh and Baishya \cite{SHAIKH2} studied
locally $\phi$-symmetric LP-Sasakian manifolds. The notion of locally $\phi$-symmetric manifolds in
different structures has been studied by several authors (see, \cite{DE}, \cite{SHAIKH2},
\cite{SABINA1}, \cite{HUI1}, \cite{TAKA}). An LP-Sasakian manifold is said to be $\phi$-symmetric \cite{DE}
if it satisfies the condition
\begin{equation}
\label{eqn1.1}
\phi^2((\nabla_W R) (X,Y)Z) = 0
\end{equation}
for arbitrary vector fields $X$, $Y$, $Z$ and $W$ on $M$.\\
\indent In particular, if $X$, $Y$, $Z$, $W$ are horizontal vector fields, i.e., orthogonal to $\xi$,
then it is called locally $\phi$-symmetric LP-Sasakian manifold \cite{TAKA}.\\
\indent It is easy to check that an LP-Sasakian manifold is $\phi$-symmetric if and only if it is
locally symmetric or locally $\phi$-symmetric.\\
\indent Recently De and Sarkar \cite{DSAR} studied $\phi$-Ricci symmetric Sasakian manifolds.
In this conncetion Shukla and Shukla \cite{SHUKLA} studied $\phi$-Ricci symmetric Kenmotsu manifolds.
An LP-Sasakian manifold is said to be $\phi$-Ricci symmetric \cite{DSAR} if it satisfies
\begin{equation}
\label{eqn1.2}
\phi^2((\nabla_X Q)(Y)) = 0,
\end{equation}
where $Q$ is the Ricci operator, i.e., $g(QX,Y) = S(X,Y)$ for all vector fields $X$, $Y$.\\
\indent If $X$, $Y$ are horizontal vector fields then the manifold is said to be locally $\phi$-Ricci symmetric.\\
\indent It is easy to check that an LP-Sasakian manifold is $\phi$-Ricci symmetric if and only if it is Ricci symmetric
or locally $\phi$-Ricci symmetric.\\
\indent The object of the present paper is to study the locally $\phi$-symmetric and locally $\phi$-Ricci symmetric
LP-Sasakian manifolds admitting semi-symmetric metric connection. The paper is organized as follows.
Section~2 is concerned with some preliminaries about LP-Sasakian manifolds and semi-symmetric metric connection.
Section~3 is devoted to the study of locally $\phi$-symmetric LP-Sasakian manifolds admitting semi-symmetric metric
connection and obtained a necessary and sufficient condition for a locally $\phi$-symmetric LP-Sasakian
manifold with respect to semi-symmetric metric connection to be locally $\phi$-symmetric LP-Sasakian
manifold with respect to Levi-Civita connection. Section 4 deals with the study of locally $\phi$-Ricci
symmetric LP-Sasakian manifolds admitting semi-symmetric metric connection.
\section{Preliminaries}
An $n$-dimensional smooth manifold $M$ is said to be an LP-Sasakian manifold (\cite{MIHAI2}, \cite{SHAIKH1})
if it admits a (1, 1) tensor field $\phi$, a unit timelike  vector field $\xi$, an 1-form $\eta$
and a Lorentzian metric $g$, which satisfy
\begin{equation}
\label{eqn2.1}
\eta(\xi)=-1,\ \ g(X,\xi)=\eta(X), \ \ \phi^{2}X=X+\eta(X)\xi,
\end{equation}
\begin{equation}
\label{eqn2.2}
g(\phi X, \phi Y) = g(X,Y)+\eta(X)\eta(Y), \ \ \nabla_{X}\xi = \phi X,
\end{equation}
\begin{equation}
\label{eqn2.3}
(\nabla_{X}\phi)(Y) = g(X,Y)\xi+\eta(Y)X+2\eta(X)\eta(Y)\xi,
\end{equation}
where $\nabla$ denotes the operator of covariant differentiation with respect to the Lorentzian metric $g$.
It can be easily seen that in an LP-Sasakian manifold, the following relations hold:
\begin{equation}
\label{eqn2.4}
\phi\xi = 0, \ \ \eta \circ \phi = 0, \ \ \ \text{rank }\phi = n-1.
\end{equation}
Again, if we take
\begin{equation*}
\Omega(X,Y) = g(X,\phi Y)
\end{equation*}
for any vector fields $X$, $Y$, then the tensor field $\Omega(X,Y)$ is a symmetric (0,2) tensor field \cite{MAT}.
Also, since the vector field $\eta$ is closed in an LP-Sasakian manifold, we have (\cite{DE0}, \cite{MAT})
\begin{equation}
\label{eqn2.5}
(\nabla_{X}\eta)(Y) = \Omega(X,Y),\ \ \ \Omega(X,\xi) = 0
\end{equation}
for any vector fields $X$ and $Y$.\\
\indent Let $M$ be an $n$-dimensional LP-Sasakian manifold with structure $(\phi,\xi,\eta,g)$. Then
the following relations hold (\cite{SHAIKH2}, \cite{SHAIKH1}):
\begin{equation}
\label{eqn2.6}
R(X,Y)\xi=\eta(Y)X-\eta(X)Y,
\end{equation}
\begin{equation}
\label{eqn2.6a}
\eta(R(X,Y)Z)=\eta(X)g(Y,Z)-\eta(Y)g(X,Z),
\end{equation}
\begin{equation}
\label{eqn2.7}
S(X,\xi) = (n-1)\eta(X),
\end{equation}
\begin{equation}
\label{eqn2.8}
S(\phi X, \phi Y) = S(X,Y) + (n-1) \eta(X)\eta(Y),
\end{equation}
\begin{equation}
\label{eqn2.8a}
(\nabla_WR)(X,Y)\xi = \Omega(Y,W)X - \Omega(X,W)Y - R(X,Y)\phi W,
\end{equation}
\begin{equation}
\label{eqn2.8b}
(\nabla_WR)(X,\xi)Y = \Omega(W,Z)X - g(X,Z)\phi W - R(X,\phi W)Z
\end{equation}
for any vector fields $X$, $Y$, $Z$, where $R$ is the curvature tensor of $g$.\\
\indent Let $M$ be an $n$-dimensional LP-Sasakian manifold and $\nabla$ be the Levi-Civita
connection on $M$. A linear connection $\widetilde{\nabla}$ on $M$ is said to be semi-symmetric
if the torsion tensor $\tau$ of the connection $\widetilde{\nabla}$
\begin{equation*}
\tau(X,Y) = \widetilde{\nabla}_{X}Y - \widetilde{\nabla}_{Y}X - [X,Y]
\end{equation*}
satisfies
\begin{equation}
\label{eqn2.9}
\tau(X,Y) = \eta(Y)X - \eta(X)Y
\end{equation}
for all $X$, $Y\in \chi(M)$; $\chi(M)$ being the Lie algebra of all smooth vector fields on $M$.
A semi-symmetric connection $\widetilde{\nabla}$ is called semi-symmetric
metric connection if it further satisfies
\begin{equation}
\label{eqn2.10}
\widetilde{\nabla}g = 0.
\end{equation}
\indent A semi-symmetric metric connection $\widetilde{\nabla}$ in an LP-Sasakian manifold
is defined by (\cite{sh},\cite{YANO5}):
\begin{equation}
\label{eqn2.11}
\widetilde{\nabla}_{X}Y = \nabla_{X}Y + \eta(Y)X - g(X,Y)\xi.
\end{equation}
\indent If $R$ and $\widetilde{R}$ are respectively the curvature tensor of the Levi-Civita connection
$\nabla$ and the semi-symmetric metric connection $\widetilde{\nabla}$ in an LP-Sasakian manifold, then we have \cite{PER}
\begin{eqnarray}
\label{eqn2.12}
\widetilde{R}(X,Y)Z &=& R(X,Y)Z - \alpha(Y,Z)X + \alpha(X,Z)Y\\
\nonumber&-&g(Y,Z)LX + g(X,Z)LY,
\end{eqnarray}
where $\alpha$ is a symmetric (0,2) tensor field given by
\begin{equation}
\label{eqn2.13}
\alpha(X,Y) = (\widetilde{\nabla}_{X}\eta)(Y) + \frac{1}{2} g(X,Y),
\end{equation}
\begin{equation}
\label{eqn2.14}
LX = \widetilde{\nabla}_{X}\xi + \frac{1}{2} X = \phi X - \frac{1}{2}X - \eta(X)\xi
\end{equation}
and
\begin{equation}
\label{eqn2.15}
g(LX,Y) = \alpha(X,Y).
\end{equation}
\begin{lem}
\emph{\cite{PER}} In an LP-Sasakian manifold with semi-symmetric metric connection $\widetilde{\nabla}$,
we have
\begin{equation}
\label{eqn2.16}
\widetilde{R}(X,Y)Z + \widetilde{R}(Y,Z)X + \widetilde{R}(Z,X)Y = 0,
\end{equation}
\begin{equation}
\label{eqn2.17}
g(\widetilde{R}(X,Y)Z,U) = - g(\widetilde{R}(Y,X)Z,U),
\end{equation}
\begin{equation}
\label{eqn2.18}
g(\widetilde{R}(X,Y)Z,U) = - g(\widetilde{R}(X,Y)U,Z),
\end{equation}
\begin{equation}
\label{eqn2.19}
g(\widetilde{R}(X,Y)Z,U) = g(\widetilde{R}(Z,U)X,Y).
\end{equation}
\end{lem}
\begin{lem}
\emph{\cite{PER}} In an $n$-dimensional LP-Sasakian manifold the Ricci tensor $\widetilde{S}$
and scalar curvature $\widetilde{r}$ with respect to semi-symmetric metric connection $\widetilde{\nabla}$ are
given by
\begin{equation}
\label{eqn2.20}
\widetilde{S}(X,Y) = S(X,Y) - (n-2)\alpha(X,Y) - a g(X,Y)
\end{equation}
and
\begin{equation}
\label{eqn2.21}
\widetilde{r} = r - 2(n-1)a,
\end{equation}
where $a = \text{ tr. } \alpha$, $S$ and $r$ denote the Ricci tensor and scalar curvature of Levi-Civita
connection $\nabla$ respectively.
\end{lem}
\begin{lem}
\emph{\cite{PER}} Let $M$ be an $n$-dimensional LP-Sasakian manifold with semi-symmetric metric
connection $\widetilde{\nabla}$. Then we have
\begin{equation}
\label{eqn2.22}
g(\widetilde{R}(X,Y)Z,\xi) = \eta(\widetilde{R}(X,Y)Z) = (\widetilde{\nabla}_{X}\eta)(Z)\eta(Y) - (\widetilde{\nabla}_{Y}\eta)(Z)\eta(X),
\end{equation}
\begin{equation}
\label{eqn2.23}
\widetilde{R}(\xi,X)\xi = - \widetilde{\nabla}_{X}\xi = X + \eta(X)\xi - \phi X,
\end{equation}
\begin{equation}
\label{eqn2.24}
\widetilde{R}(X,Y)\xi = \eta(X)\widetilde{\nabla}_{Y}\xi - \eta(Y)\widetilde{\nabla}_{X}\xi,
\end{equation}
\begin{equation}
\label{eqn2.25}
\widetilde{R}(\xi,X)Y = \eta(Y)\widetilde{\nabla}_{X}\xi - g(Y,\widetilde{\nabla}_{X}\xi)\xi,
\end{equation}
\begin{equation}
\label{eqn2.26}
\widetilde{S}(X,\xi) = \Big(\frac{n}{2}-a\Big)\eta(X),
\end{equation}
\begin{eqnarray}
\label{eqn2.27}
\widetilde{S}(\phi X,\phi Y) &=& S(X,Y) + \Big(\frac{n}{2}-a\Big)\eta(X)\eta(Y)\\
\nonumber&-&(n-2)\alpha(X,Y) - a g(X,Y)
\end{eqnarray}
for arbitrary vector fields $X$, $Y$ and $Z$.
\end{lem}
From (\ref{eqn2.2}), (\ref{eqn2.3}), (\ref{eqn2.5}), (\ref{eqn2.11}) and (\ref{eqn2.14}), we get
\begin{eqnarray}
\label{eqn2.28}
(\widetilde{\nabla}_{W} \widetilde{R})(X,Y)\xi &=& R(X,Y)W - R(X,Y)\phi W + \alpha(X,W)Y\\
\nonumber&-& \alpha(Y,W)X + g(X,W)LY - g(Y,W)LX\\
\nonumber&+& \alpha(Y,\phi W)X - \alpha(X,\phi W)Y + \Omega(Y,W)LX\\
\nonumber&-& \Omega(X,W)LY + g(X,W)Y - g(Y,W)X\\
\nonumber&+& g(Y,W)\phi X - g(X,W)\phi Y + \Omega(Y,W)X \\
\nonumber&-& \Omega(X,W)Y + \Omega(X,W)\phi Y - \Omega(Y,W)\phi X\\
\nonumber&+& \eta(X)[g(Y,W)-\Omega(Y,W)]\xi\\
\nonumber&-& \eta(Y)[g(X,W)-\Omega(X,W)]\xi
\end{eqnarray}
for arbitrary vector fields $X$, $Y$ and $W$. Also from (\ref{eqn2.11}), (\ref{eqn2.12}) and (\ref{eqn2.18}), we have
\begin{equation}
\label{eqn2.29}
g((\widetilde{\nabla}_{W} \widetilde{R})(X,Y)Z,U) = - g((\widetilde{\nabla}_{W} \widetilde{R})(X,Y)U,Z).
\end{equation}
From (\ref{eqn2.14}) we have
\begin{equation}
\label{eqn2.33}
\alpha(X,\xi) = \frac{1}{2}\eta(X),
\end{equation}
\begin{equation}
\label{eqn2.34}
(\nabla_{W}\alpha)(X,\xi) = \frac{1}{2}\Omega(W,X) - \alpha(X,\phi W),
\end{equation}
\begin{eqnarray}
\label{eqn2.35}
(\nabla_{W}L)(X) &=& [g(W,X) - \Omega(W,X)]\xi\\
\nonumber&+&\eta(X)[W - \phi W] + 2\eta(X)\eta(W)\xi.
\end{eqnarray}
Again by virtue of (\ref{eqn2.33}) - (\ref{eqn2.35}) we have from (\ref{eqn2.11}) and (\ref{eqn2.12}) that
\begin{eqnarray}
\label{eqn2.30}
&&(\widetilde{\nabla}_{W} \widetilde{R})(X,Y)Z\\
\nonumber&=& (\nabla_{W}R)(X,Y)Z - g(R(X,Y)Z,W)\xi + [g(W,Y) - \Omega(W,Y)]\eta(Z)X\\
\nonumber&+& [g(W,Z) - \Omega(W,Z)]\eta(Y)X + 2\eta(Z)\eta(W)[\eta(Y)X-\eta(X)Y]\\
\nonumber&+&\alpha(Y,Z)[g(X,W)\xi - \eta(X)W] + [\Omega(W,X) - g(W,X)]\eta(Z)Y\\
\nonumber&+& [\Omega(W,Z) - g(W,Z)]\eta(X)Y + \alpha(X,Z)[\eta(Y)W - g(Y,W)\xi]\\
\nonumber&-& g(Y,Z)[\{g(W,X) - \Omega(W,X) - \alpha(X,W)\}\xi +\eta(X)\{\frac{1}{2}W - \phi W + 2 \eta(W)\xi\}]\\
\nonumber&+& g(X,Z)[\{g(W,Y) - \Omega(W,Y) - \alpha(Y,W)\}\xi +\eta(Y)\{\frac{1}{2}W - \phi W + 2 \eta(W)\xi\}].
\end{eqnarray}
By virtue of (\ref{eqn2.33}) and (\ref{eqn2.35}) it follows from
(\ref{eqn2.11}) that
\begin{eqnarray}
\label{eqn2.37} (\widetilde{\nabla}_{X}\widetilde{S})(Y,Z)
&=&(\nabla_{X}S)(Y,Z) - [S(X,Y)+\alpha(X,Y)]\eta(Z)\\
\nonumber&+&[\frac{3}{2}g(X,Z)+(n-2)\Omega(X,Z)]\eta(Y)\\
\nonumber&-&(n-2)[g(X,Y)-\Omega(X,Y)]\eta(Z) - da(X)g(Y,Z).
\end{eqnarray}
Also from (\ref{eqn2.7}) we have
\begin{equation}
\label{eqn2.38} (\nabla_{X}S)(Y,\xi) = (n-1)\Omega(X,Y) - S(Y,\phi
X).
\end{equation}
\section{Locally $\phi$-symmetric LP-Sasakian manifolds admitting semi-symmetric metric connection}
\begin{df}
An LP-Sasakian manifold $M$ is said to be locally $\phi$-symmetric with respect to semi-symmetric metric
connection if its curvature tensor $\widetilde{R}$ satisfies the condition
\begin{equation}
\label{eqn3.1}
\phi^2((\widetilde{\nabla}_W \widetilde{R}) (X,Y)Z) = 0
\end{equation}
for all horizontal vector fields $X$, $Y$, $Z$ and $W$.
\end{df}
We now consider a locally $\phi$-symmetric LP-Sasakian manifold with respect to semi-symmetric metric connection.
Then by virtue of (\ref{eqn2.1}) it follows from (\ref{eqn3.1}) that
\begin{equation}
\label{eqn3.2}
(\widetilde{\nabla}_W \widetilde{R}) (X,Y)Z + \eta((\widetilde{\nabla}_W \widetilde{R}) (X,Y)Z)\xi = 0.
\end{equation}
Using (\ref{eqn2.29}) in (\ref{eqn3.2}), we get
\begin{equation}
\label{eqn3.3}
(\widetilde{\nabla}_W \widetilde{R}) (X,Y)Z = g((\widetilde{\nabla}_W \widetilde{R}) (X,Y)\xi, Z)\xi.
\end{equation}
In view of (\ref{eqn2.28}) it follows from (\ref{eqn3.3}) that
\begin{eqnarray}
\label{eqn3.4}
\nonumber (\widetilde{\nabla}_{W} \widetilde{R})(X,Y)Z &=& \Big[g(R(X,Y)W,Z) - g(R(X,Y)\phi W,Z) + \alpha(X,W)g(Y,Z)\\
&-& \alpha(Y,W)g(X,Z) + g(X,W)\alpha(Y,Z) - g(Y,W)\alpha(X,Z)\\
\nonumber&+& \alpha(Y,\phi W)g(X,Z)- \alpha(X,\phi W)g(Y,Z) + \Omega(Y,W)\alpha(X,Z)\\
\nonumber&-& \Omega(X,W)\alpha(Y,Z) + g(X,W)g(Y,Z) - g(Y,W)g(X,Z)\\
\nonumber&+& g(Y,W)\Omega(X,Z) - g(X,W)\Omega(Y,Z) + \Omega(Y,W)g(X,Z)\\
\nonumber&-& \Omega(X,W)g(Y,Z) + \Omega(X,W)\Omega(Y,Z) - \Omega(Y,W)\Omega(X,Z)\Big]\xi
\end{eqnarray}
for all horizontal vector fields $X$, $Y$, $Z$ and $W$. Next, let us assume that in an LP-Sasakian manifold, the relation (\ref{eqn3.4}) holds
for all horizontal vector fields $X$, $Y$, $Z$ and $W$. Then it follows from (\ref{eqn2.30}) that (\ref{eqn3.4}) holds
and consequently the manifold is locally $\phi$-symmetric with respect to semi-symmetric metric connection. This leads to the following:
\begin{thm}
An LP-Sasakian manifold is locally $\phi$-symmetric with respect to semi-symmetric metric connection if
and only if the relation \emph{(\ref{eqn3.4})} holds for all horizontal vector fields $X$, $Y$, $Z$ and $W$.
\end{thm}
In view of (\ref{eqn2.29}), it follows from (\ref{eqn3.2}) that
\begin{equation}
\label{eqn3.5}
(\widetilde{\nabla}_W \widetilde{R}) (X,Y)\xi = 0.
\end{equation}
From (\ref{eqn2.28}) and (\ref{eqn3.5}) it follows that
\begin{eqnarray}
\label{eqn3.6}
&&R(X,Y)W - R(X,Y)\phi W\\
\nonumber&=& g(Y,W)X - g(X,W)Y + g(X,W)\phi Y - g(Y,W)\phi X\\
\nonumber&+& \Omega(X,W)Y - \Omega(Y,W)X + \Omega(Y,W)\phi X - \Omega(X,W) \phi Y\\
\nonumber&+& \alpha(Y,W)X - \alpha(X,W)Y + g(Y,W)LX - g(X,W)LY\\
\nonumber&+& \alpha(X,\phi W)Y- \alpha(Y,\phi W)X + \Omega(X,W)LY - \Omega(Y,W)LX
\end{eqnarray}
for horizontal vector fields $X$, $Y$ and $W$. Contracting (\ref{eqn3.6}), we get
\begin{eqnarray}
\label{eqn3.7}
S(Y,W) - S(Y,\phi W)&=& (n-1+a-\psi)[g(Y,W) - \Omega(Y,W)]\\
\nonumber&+& (n-2)[\alpha(Y,W) - \alpha(Y,\phi W)],
\end{eqnarray}
where $\psi = \text{ tr. } \Omega$ and $a = \text{ tr. } \alpha$. Hence we can state the following:
\begin{thm}
In a locally $\phi$-symmetric LP-Sasakian manifold with semi-symmetric metric connection the curvature tensor and the Ricci tensor are respectively given by \emph{(\ref{eqn3.6})} and \emph{(\ref{eqn3.7})}.
\end{thm}
We now consider a locally $\phi$-symmetric LP-Sasakian manifold with Levi-Civita connection. Then in
\cite{SHAIKH2}, Shaikh and Baishya proved that
\begin{thm}
An LP-Sasakian manifold $(M^n,g)$ is locally $\phi$-symmetric with respect to
Levi-Civita connection if and only if the following relation
\begin{eqnarray}
\label{eqn3.12}
&&(\nabla_{W}R)(X,Y)Z\\
\nonumber&=&\big[2\{\Omega(Y,W)g(X,Z) - \Omega(X,W)g(Y,Z)\} + \Omega(Y,Z)g(X,W) - \Omega(X,Z)g(Y,W)\\
\nonumber&+&2\{\Omega(Y,Z)\eta(X)\eta(W) - \Omega(X,Z)\eta(Y)\eta(W)\} - g(\phi R(X,Y)W,Z)\big]\xi\\
\nonumber&+&\eta(X)[\Omega(W,Z)Y - g(Y,Z)\phi W - R(Y,\phi W)Z] - \eta(Y)[\Omega(W,Z)X - g(X,Z)\phi W\\
\nonumber&-&R(X,\phi W)Z] - \eta(Z)[2\{\Omega(Y,W)X - \Omega(X,W)Y\} - \phi R(X,Y)W - g(Y,W)\phi X\\
\nonumber&+&g(X,W) \phi Y] + 2\{\eta(Y)\phi X - \eta(X)\phi Y\}\eta(Z)\eta(W).
\end{eqnarray}
\end{thm}
holds for arbitrary vector fields $X$, $Y$, $Z$, $W \in \chi(M)$.\\
\indent Now we take a locally $\phi$-symmetric LP-Sasakian manifold with respect to semi-symmetric metric
connection. Then the relation (\ref{eqn3.4}) holds for any horizontal vector fields $X$, $Y$, $Z$, $W$.\\
\indent Let $X$, $Y$, $Z$, $W$ be arbitrary vector fields of $\chi(M)$. We now compute\\
$(\widetilde{\nabla}_{\phi^2W} \widetilde{R})(\phi^2 X, \phi^2 Y)\phi^2 Z$ in two different ways. Firstly by virtue of (\ref{eqn2.1})
it follows from (\ref{eqn3.4}) that
\begin{eqnarray}
\label{eqn3.13}
&&(\widetilde{\nabla}_{\phi^2W} \widetilde{R})(\phi^2 X, \phi^2 Y)\phi^2 Z\\
\nonumber&=&\big[g(R(\phi^2X,\phi^2Y)\phi^2W,\phi^2Z) - g(R(\phi^2X,\phi^2Y)\phi^3W,\phi^2Z)\\
\nonumber&+&\alpha(\phi^2X,\phi^2W)\{g(Y,Z)+\eta(Y)\eta(Z)\}-\alpha(\phi^2Y,\phi^2W)\{g(X,Z)+\eta(X)\eta(Z)\}\\
\nonumber&+&\alpha(\phi^2Y,\phi^2Z)\{g(X,W)+\eta(X)\eta(W)\} - \alpha(\phi^2X,\phi^2Z)\{g(Y,W)+\eta(Y)\eta(W)\}\\
\nonumber&+&\alpha(\phi^2Y,\phi^3W)\{g(X,Z)+\eta(X)\eta(Z)\} - \alpha(\phi^2X,\phi^3W)\{g(Y,Z)+\eta(Y)\eta(Z)\}\\
\nonumber&+&\Omega(Y,W)\alpha(\phi^2X,\phi^2Z) - \Omega(X,W)\alpha(\phi^2Y,\phi^2Z)+\{g(X,W)+\eta(X)\eta(W)\}\{g(Y,Z)\\
\nonumber&+&\eta(Y)\eta(Z)\}-\{g(Y,W)+\eta(Y)\eta(W)\}\{g(X,Z)+\eta(X)\eta(Z)\}+\{g(Y,W)\\
\nonumber&+&\eta(Y)\eta(W)\}\Omega(X,Z)-\{g(X,W)+\eta(X)\eta(W)\}\Omega(Y,Z)\\
\nonumber&+&\{g(X,Z)+\eta(X)\eta(Z)\}\Omega(Y,W)-\{g(Y,Z)+\eta(Y)\eta(Z)\}\Omega(X,W)\\
\nonumber&+&\Omega(X,W)\Omega(Y,Z) - \Omega(Y,W)\Omega(X,Z)\big]\xi.
\end{eqnarray}
From (\ref{eqn2.4}) we have
\begin{equation}
\label{eqn3.14} g(\phi^2X,\xi)=g(\phi^2Y,\xi)=g(\phi^2Z,\xi)=0
\end{equation}
and hence $\phi^2X$, $\phi^2Y$, $\phi^2Z$ are horizontal vector fields of $\chi(M)$. Then by virtue of (\ref{eqn2.1})
it follows that
\begin{eqnarray}
\label{eqn3.15}
R(\phi^2X,\phi^2Y)\phi^2W &=&R(X,Y)W + \{\eta(Y)X - \eta(X)Y\}\eta(W)\\
\nonumber&+&\{g(Y,W)\eta(X) - g(X,W)\eta(Y)\}\xi,
\end{eqnarray}
\begin{equation}
\label{eqn3.16}
R(\phi^2X,\phi^2Y)\phi^3W = R(X,Y)\phi W +\{\Omega(Y,W)\eta(X) - \Omega(X,W)\eta(Y)\}\xi,
\end{equation}
\begin{equation}
\label{eqn3.17}
\alpha(\phi^2X, \phi^2W) = \alpha(X,W) + \frac{1}{2}\eta(X)\eta(W).
\end{equation}
In view of (\ref{eqn3.15}) - (\ref{eqn3.17}), (\ref{eqn3.13}) yields
\begin{eqnarray}
\label{eqn3.20}
&&\ \ \ \ (\widetilde{\nabla}_{\phi^2W} \widetilde{R})(\phi^2 X, \phi^2 Y)\phi^2 Z\\
\nonumber&=&\big[g(R(X,Y)W,Z) - g(R(X,Y)\phi W,Z)+\alpha(X,W)\{g(Y,Z)+\eta(Y)\eta(Z)\}\\
\nonumber&-&\alpha(Y,W)\{g(X,Z)+\eta(X)\eta(Z)\}+\frac{1}{2}\{\eta(X)g(Y,Z)-\eta(Y)g(X,Z)\}\eta(W)\\
\nonumber&+&\alpha(Y,Z)g(X,W)-\alpha(X,Z)g(Y,W)+\frac{1}{2}\{\eta(Y)g(X,W) - \eta(X)g(Y,W)\}\eta(Z)\\
\nonumber&+&\{\eta(X)\alpha(Y,Z)-\eta(Y)\alpha(X,Z)\}\eta(W) + \alpha(Y,\phi W)g(X,Z)-\alpha(X,\phi W)g(Y,Z)\\
\nonumber&+&\{\eta(X)\alpha(Y,\phi W) - \eta(Y)\alpha(X,\phi W)\}\eta(Z)+\Omega(Y,W)\alpha(X,Z) - \Omega(X,W)\alpha(Y,Z)\\
\nonumber&+&\frac{1}{2}\{\eta(X)\Omega(Y,W)-\eta(Y)\Omega(X,W)\}\eta(Z)+g(X,W)g(Y,Z) - g(Y,W)g(X,Z)\\
\nonumber&+&g(Y,W)\Omega(X,Z)-g(X,W)\Omega(Y,Z)+\{\eta(Y)\Omega(X,Z) - \eta(X)\Omega(Y,Z)\}\eta(W)\\
\nonumber&+&\Omega(Y,W)g(X,Z)-\Omega(X,W)g(Y,Z)+\{\eta(X)\Omega(Y,W)-\eta(Y)\Omega(X,W)\}\eta(Z)\\
\nonumber&+&\Omega(X,W)\Omega(Y,Z) - \Omega(Y,W)\Omega(X,Z)\big]\xi.
\end{eqnarray}
By virtue of (\ref{eqn2.1}) we have
\begin{eqnarray}
\label{eqn3.21} (\widetilde{\nabla}_{\phi^2W}
\widetilde{R})(\phi^2 X, \phi^2 Y)\phi^2 Z &=&
(\widetilde{\nabla}_{W} \widetilde{R})(\phi^2 X, \phi^2 Y)\phi^2
Z\\
\nonumber&+& \eta(W)(\widetilde{\nabla}_{\xi}
\widetilde{R})(\phi^2 X, \phi^2 Y)\phi^2 Z.
\end{eqnarray}
Now for any horizontal vector fields $X$, $Y$ and $Z$ we have from
(\ref{eqn3.4}) that
\begin{equation}
\label{eqn3.22} (\widetilde{\nabla}_{\xi}\widetilde{R})(X,Y)Z = 0,
\end{equation}
which implies that
\begin{equation}
\label{eqn3.23} (\widetilde{\nabla}_{\xi}\widetilde{R})(\phi^2 X,
\phi^2 Y)\phi^2 Z = 0.
\end{equation}
Using (\ref{eqn3.23}) in (\ref{eqn3.21}) we obtain
\begin{equation}
\label{eqn3.24} (\widetilde{\nabla}_{\phi^2W}
\widetilde{R})(\phi^2 X, \phi^2 Y)\phi^2 Z =
(\widetilde{\nabla}_{W} \widetilde{R})(\phi^2 X, \phi^2 Y)\phi^2
Z.
\end{equation}
In view of (\ref{eqn2.1}), we have
\begin{eqnarray}
\label{eqn3.25} &&(\widetilde{\nabla}_{W}\widetilde{R})(\phi^2X,\phi^2Y)\phi^2Z\\
\nonumber&=&(\widetilde{\nabla}_{W}\widetilde{R})(X,Y)Z + \eta(Z)(\widetilde{\nabla}_{W}\widetilde{R})(X,Y)\xi\\
\nonumber&+&\eta(Y)(\widetilde{\nabla}_{W}\widetilde{R})(X,\xi)Z + \eta(Y)\eta(Z)(\widetilde{\nabla}_{W}\widetilde{R})(X,\xi)\xi\\
\nonumber&+&\eta(X)(\widetilde{\nabla}_{W}\widetilde{R})(\xi,Y)Z + \eta(X)\eta(Z)(\widetilde{\nabla}_{W}\widetilde{R})(\xi,Y)Z.
\end{eqnarray}
Using (\ref{eqn2.30}) in (\ref{eqn3.25}) we get
\begin{eqnarray}
\label{eqn3.31} &&(\widetilde{\nabla}_{W}\widetilde{R})(\phi^2X,\phi^2Y)\phi^2Z\\
\nonumber&=&(\widetilde{\nabla}_{W}\widetilde{R})(X,Y)Z - \eta(Z)R(X,Y)\phi W - \eta(Y)R(X,\phi W)Z + \eta(X)R(Y,\phi W)Z\\
\nonumber&+&\frac{1}{2}\big[\eta(Z)\{\Omega(Y,W)X-\Omega(X,W)Y\}+\eta(Y)\Omega(W,Z)X - \eta(X)\Omega(W,Z)Y\big]\\
\nonumber&-&\eta(Z)\{\alpha(Y,\phi W)X - \alpha(X,\phi W)Y\}+\eta(Y)\alpha(Z,\phi W)X - \eta(X)\alpha(Z,\phi W)Y\\
\nonumber&-&\eta(X)\{\alpha(Y,Z)W - \eta(W)\alpha(Y,Z)\xi\}+\eta(Y)\eta(Z)\{\alpha(X,W)\xi - \alpha(X,\phi W)\xi\}\\
\nonumber&-&\eta(X)\eta(Z)\{\alpha(Y,W)\xi - \alpha(Y,\phi W)\xi\} - \frac{1}{2}\eta(X)\{g(Y,Z)W - \eta(W)g(Y,Z)\xi\}\\
\nonumber&+&\frac{1}{2}\eta(Y)\{g(X,Z)W - \eta(W)g(X,Z)\xi\}.
\end{eqnarray}
From (\ref{eqn3.24}) and (\ref{eqn3.31}), we get
\begin{eqnarray}
\label{eqn3.32} &&(\widetilde{\nabla}_{\phi^2W}\widetilde{R})(\phi^2X,\phi^2Y)\phi^2Z\\
\nonumber&=&(\widetilde{\nabla}_{W}\widetilde{R})(X,Y)Z - \eta(Z)R(X,Y)\phi W - \eta(Y)R(X,\phi W)Z + \eta(X)R(Y,\phi W)Z\\
\nonumber&+&\frac{1}{2}\big[\eta(Z)\{\Omega(Y,W)X-\Omega(X,W)Y\}+\eta(Y)\Omega(W,Z)X - \eta(X)\Omega(W,Z)Y\big]\\
\nonumber&-&\eta(Z)\{\alpha(Y,\phi W)X - \alpha(X,\phi W)Y\}+\eta(Y)\alpha(Z,\phi W)X - \eta(X)\alpha(Z,\phi W)Y\\
\nonumber&-&\eta(X)\{\alpha(Y,Z)W - \eta(W)\alpha(Y,Z)\xi\}+\eta(Y)\eta(Z)\{\alpha(X,W)\xi - \alpha(X,\phi W)\xi\}\\
\nonumber&-&\eta(X)\eta(Z)\{\alpha(Y,W)\xi - \alpha(Y,\phi W)\xi\} - \frac{1}{2}\eta(X)\{g(Y,Z)W - \eta(W)g(Y,Z)\xi\}\\
\nonumber&+&\frac{1}{2}\eta(Y)\{g(X,Z)W - \eta(W)g(X,Z)\xi\}.
\end{eqnarray}
From (\ref{eqn3.20}) and (\ref{eqn3.32}) we obtain
\begin{eqnarray}
\label{eqn3.33} &&(\widetilde{\nabla}_{W}\widetilde{R})(X,Y)Z\\
\nonumber&=&\Big[g(R(X,Y)W,Z) - g(R(X,Y)\phi W,Z) + \alpha(X,W)g(Y,Z) - \alpha(Y,W)g(X,Z)\\
\nonumber&+&\alpha(Y,Z)g(X,W)-\alpha(X,Z)g(Y,W)+\frac{1}{2}\{\eta(Y)g(X,W) - \eta(X)g(Y,W)\}\eta(Z)\\
\nonumber&-&\eta(Y)\eta(W)\alpha(X,Z)+\alpha(Y,\phi W)g(X,Z)-\alpha(X,\phi W)g(Y,Z)\\
\nonumber&+&\Omega(Y,W)\alpha(X,Z) - \Omega(X,W)\alpha(Y,Z)+\frac{1}{2}\{\eta(X)\Omega(Y,W)-\eta(Y)\Omega(X,W)\}\eta(Z)\\
\nonumber&+&g(X,W)g(Y,Z) - g(Y,W)g(X,Z)+g(Y,W)\Omega(X,Z)-g(X,W)\Omega(Y,Z)\\
\nonumber&+&\{\eta(Y)\Omega(X,Z)-\eta(X)\Omega(Y,Z)\}\eta(W) + \Omega(Y,W)g(X,Z)-\Omega(X,W)g(Y,Z)
\end{eqnarray}
\begin{eqnarray}
\nonumber&+&\{\eta(X)\Omega(Y,W) - \eta(Y)\Omega(X,W)\}\eta(Z)+\Omega(X,W)\Omega(Y,Z)-\Omega(Y,W)\Omega(X,Z)\Big]\xi\\
\nonumber&+&\eta(Z)R(X,Y)\phi W + \eta(Y)R(X,\phi W)Z -\eta(X)R(Y,\phi W)Z\\
\nonumber&-&\frac{1}{2}\big[\eta(Z)\{\Omega(Y,W)X - \Omega(X,W)Y\} + \eta(Y)\Omega(W,Z)X - \eta(X)\Omega(W,Z)Y\big]\\
\nonumber&+&\eta(Z)\{\alpha(Y,\phi W)X - \alpha(X,\phi W)Y\} - \eta(Y)\alpha(Z,\phi W)X + \eta(X)\alpha(Z,\phi W)Y\\
\nonumber&+&\eta(X)\alpha(Y,Z)W+\frac{1}{2}\{\eta(X)g(Y,Z)W - \eta(Y)g(X,Z)W\}.
\end{eqnarray}
\indent Thus in a locally $\phi$-symmetric LP-Sasakian manifold with respect to semi-symmetric metric connection, the relation
(\ref{eqn3.33}) holds for any $X$, $Y$, $Z$, $W\in\chi(M)$.\\
\indent Next, if the relation (\ref{eqn3.33}) holds in an LP-Sasakian manifold with respect to semi-symmetric metric connection
then for any horizontal vector fields $X$, $Y$, $Z$, $W$, we obtain the relation (\ref{eqn3.4}) and hence the manifold is
locally $\phi$-symmetric with respect to semi-symmetric metric connection. Thus we can state the following:
\begin{thm}
An LP-Sasakian manifold $(M^n,g)$ is locally $\phi$-symmetric with respect to semi-symmetric metric connection
if and only if the relation \emph{(\ref{eqn3.33})} holds for any vector fields $X$, $Y$, $Z$, $W\in\chi(M)$.
\end{thm}
In view of (\ref{eqn2.30}), (\ref{eqn3.33}) yields
\begin{eqnarray}
\label{eqn3.34} &&(\nabla_{W}R)(X,Y)Z\\
\nonumber&=&[\Omega(W,Y)-g(W,Y)]\eta(Z)X+[\Omega(W,Z)-g(W,Z)]\eta(Y)X\\
\nonumber&+&2\eta(Z)\eta(W)[\eta(X)Y-\eta(Y)X]+[\alpha(Y,Z)\eta(X) - \alpha(X,Z)\eta(Y)]W\\
\nonumber&+&[g(W,X)-\Omega(W,X)]\eta(Z)Y + [g(W,Z) - \Omega(W,Z)]\eta(X)Y\\
\nonumber&+&g(Y,Z)\eta(X)[\frac{1}{2}W-\phi W]-g(X,Z)\eta(Y)[\frac{1}{2}W-\phi W]\\
\nonumber&+&\eta(Z)R(X,Y)\phi W + \eta(Y)R(X,\phi W)Z - \eta(X)R(Y,\phi W)Z\\
\nonumber&+&\frac{1}{2}\big[\eta(X)\Omega(W,Z)Y - \eta(Y)\Omega(W,Z)X - \eta(Z)\{\Omega(Y,W)X - \Omega(X,W)Y\}\big]\\
\nonumber&+&\eta(Z)\{\alpha(Y,\phi W)X - \alpha(X,\phi W)Y\}-\eta(Y)\alpha(Z,\phi W)X + \eta(X)\alpha(Z,\phi W)Y\\
\nonumber&+&\eta(X)\alpha(Y,Z)W + \frac{1}{2}\{\eta(X)g(Y,Z)W - \eta(Y)g(X,Z)W\}\\
\nonumber&+&\Big[2\{\eta(X)g(Y,Z) - \eta(Y)g(X,Z)\}\eta(Y) - g(R(X,Y)\phi W,Z)\\
\nonumber&+&\frac{1}{2}\{\eta(Y)g(X,W) - \eta(X)g(Y,W)\}\eta(Z)\\
\nonumber&-&\eta(Y)\eta(W)\alpha(X,Z)+\alpha(Y,\phi W)g(X,Z) - \alpha(X,\phi W)g(Y,Z)+\Omega(Y,W)\alpha(X,Z)
\end{eqnarray}
\begin{eqnarray}
\nonumber&-&\Omega(X,W)\alpha(Y,Z) + \frac{3}{2}\{\eta(X)\Omega(Y,W) - \eta(Y)\Omega(X,W)\}\eta(Z)+g(Y,W)\Omega(X,Z)\\
\nonumber&-&g(X,W)\Omega(Y,Z)+\{\eta(Y)\Omega(X,Z) - \eta(X)\Omega(Y,Z)\}\eta(W) + 2\{\Omega(Y,W)g(X,Z)\\
\nonumber&-&\Omega(X,W)g(Y,Z)\}+\Omega(X,W)\Omega(Y,Z) - \Omega(Y,W)\Omega(X,Z)\Big]\xi.
\end{eqnarray}
This leads to the following:
\begin{thm}
In a locally $\phi$-symmetric LP-Sasakian manifold with respect to semi-symmetric metric connection the relation
\emph{(\ref{eqn3.34})} holds for any vector fields $X$, $Y$, $Z$, $W\in\chi(M)$.
\end{thm}
From (\ref{eqn3.12}) and (\ref{eqn3.34}), we can state the following:
\begin{thm}
A locally $\phi$-symmetric LP-Sasakian manifold is invariant under a semi-symmetric metric
connection if and only if the relation
\begin{eqnarray}
\label{eqn3.35} &&[\Omega(W,Y)-g(W,Y)]\eta(Z)X+[\Omega(W,Z)-g(W,Z)]\eta(Y)X\\
\nonumber&+&2\eta(Z)\eta(W)[\eta(X)Y-\eta(Y)X]+[\alpha(Y,Z)\eta(X) - \alpha(X,Z)\eta(Y)]W\\
\nonumber&+&[g(W,X)-\Omega(W,X)]\eta(Z)Y + [g(W,Z) - \Omega(W,Z)]\eta(X)Y\\
\nonumber&+&\frac{1}{2}[g(Y,Z)\eta(X)- g(X,Z)\eta(Y)]W+\eta(Z)[g(Y,W)\phi X - g(X,W)\phi Y]\\
\nonumber&-&\frac{1}{2}\big[\eta(X)\Omega(W,Z)Y - \eta(Y)\Omega(W,Z)X - \eta(Z)\{\Omega(Y,W)X - \Omega(X,W)Y\}\big]\\
\nonumber&+&\eta(Z)\{\alpha(Y,\phi W)X - \alpha(X,\phi W)Y\}-\eta(Y)\alpha(Z,\phi W)X + \eta(X)\alpha(Z,\phi W)Y\\
\nonumber&+&\eta(X)\alpha(Y,Z)W + \frac{1}{2}\{\eta(X)g(Y,Z)W - \eta(Y)g(X,Z)W\}\\
\nonumber&+&\Big[2\{\eta(X)g(Y,Z) - \eta(Y)g(X,Z)\}\eta(Y)+\frac{1}{2}\{\eta(Y)g(X,W) - \eta(X)g(Y,W)\}\eta(Z)\\
\nonumber&-&\eta(Y)\eta(W)\alpha(X,Z)+\alpha(Y,\phi W)g(X,Z) -
\alpha(X,\phi W)g(Y,Z)+\Omega(Y,W)\alpha(X,Z)\\
\nonumber&-&\Omega(X,W)\alpha(Y,Z) + \frac{3}{2}\{\eta(X)\Omega(Y,W) - \eta(Y)\Omega(X,W)\}\eta(Z)+g(Y,W)\Omega(X,Z)\\
\nonumber&-&g(X,W)\Omega(Y,Z)+\{\eta(Y)\Omega(X,Z) - \eta(X)\Omega(Y,Z)\}\eta(W) + \Omega(Y,W)g(X,Z)\\
\nonumber&-&\Omega(X,W)g(Y,Z)+\Omega(X,W)\Omega(Y,Z) -
\Omega(Y,W)\Omega(X,Z)\Big]\xi = 0
\end{eqnarray}
holds for arbitrary vector fields $X$, $Y$, $Z$, $W\in\chi(M)$.
\end{thm}
\section{Locally $\phi$-Ricci symmetric LP-Sasakian manifolds admitting\\ semi-symmetric connection}
\begin{df}
An LP-Sasakian manifold $M$ is said to be locally $\phi$-Ricci symmetric with respect to semi-symmetric metric connection
if its satisfies the condition
\begin{equation}
\label{eqn4.1}
\phi^2((\widetilde{\nabla}_X \widetilde{Q}) (Y)) = 0
\end{equation}
for horizontal vector fields $X$ and $Y$, where $\widetilde{Q}$ is
the Ricci-operator with respect to semi-symmetric metric
connection, i.e. $g(\widetilde{Q}X,Y) = \widetilde{S}(X,Y)$ for
all vector fields $X$, $Y$.
\end{df}
Let us take an LP-Sasakian manifold, which is $\phi$-Ricci symmetric with respect to semi-symmetric metric connection.
Then by virtue of (\ref{eqn2.1}) it follows from (\ref{eqn4.1}) that
\begin{equation*}
(\widetilde{\nabla}_X \widetilde{Q}) (Y)  + \eta((\widetilde{\nabla}_X \widetilde{Q}) (Y)) \xi= 0
\end{equation*}
from which it follows that
\begin{equation}
\label{eqn4.2} (\widetilde{\nabla}_{X}\widetilde{S})(Y,Z) = 0
\end{equation}
for all horizontal vector fields $X$ and $Y$ and $Z$.\\
\indent Let $X$, $Y$, $Z$ be arbitrary vector fields of $\chi(M)$. We now compute $(\widetilde{\nabla}_{\phi^2X}\widetilde{S})(\phi^2Y, \phi^2Z)$
in two different ways. Since $\phi^2X$, $\phi^2Y$, $\phi^2Z$ are horizontal vector fields for all $X$, $Y$, $Z\in\chi(M)$, from (\ref{eqn4.2}) we have
\begin{equation}
\label{eqn4.3}
(\widetilde{\nabla}_{\phi^2X}\widetilde{S})(\phi^2Y,\phi^2Z) = 0
\end{equation}
for all $X$, $Y$, $Z\in\chi(M)$. By virtue of (\ref{eqn2.1}) we get
\begin{equation}
\label{eqn4.4}
(\widetilde{\nabla}_{\phi^2X}\widetilde{S})(\phi^2Y,\phi^2Z) = (\widetilde{\nabla}_{X}\widetilde{S})(\phi^2Y,\phi^2Z) + \eta(X)(\widetilde{\nabla}_{\xi}\widetilde{S})(\phi^2Y,\phi^2Z).
\end{equation}
Now for any horizontal vector fields $Y$ and $Z$ we have from (\ref{eqn4.2}) that
\begin{equation*}
(\widetilde{\nabla}_{\xi}\widetilde{S})(Y,Z) = 0,
\end{equation*}
which implies that
\begin{equation}
\label{eqn4.5}
(\widetilde{\nabla}_{\xi}\widetilde{S})(\phi^2Y,\phi^2Z) = 0
\end{equation}
for arbitrary vector fields $Y$, $Z\in\chi(M)$.\\
Using (\ref{eqn4.5}) in (\ref{eqn4.4}) we get
\begin{equation}
\label{eqn4.6}
(\widetilde{\nabla}_{\phi^2X}\widetilde{S})(\phi^2Y,\phi^2Z) = (\widetilde{\nabla}_{X}\widetilde{S})(\phi^2Y,\phi^2Z).
\end{equation}
In view of (\ref{eqn2.1}), we get
\begin{eqnarray}
\label{eqn4.7}
(\widetilde{\nabla}_{X}\widetilde{S})(\phi^2Y,\phi^2Z) &=& (\widetilde{\nabla}_{X}\widetilde{S})(Y,Z) + \eta(Y)(\widetilde{\nabla}_{X}\widetilde{S})(Z,\xi)\\
\nonumber&+&\eta(Z)(\widetilde{\nabla}_{X}\widetilde{S})(Z,\xi) + \eta(Y)\eta(Z)(\widetilde{\nabla}_{X}\widetilde{S})(\xi,\xi).
\end{eqnarray}
Using (\ref{eqn2.37}) in (\ref{eqn4.7}) we get
\begin{eqnarray}
\label{eqn4.8}
(\widetilde{\nabla}_{X}\widetilde{S})(\phi^2Y,\phi^2Z) &=& (\nabla_{X}S)(Y,Z) - \eta(Z)S(Y,\phi X)\\
\nonumber&+&\eta(Y)[S(X,Z) - S(Z,\phi X)]+\eta(Y)\alpha(X,Z)\\
\nonumber&+&[(2n-1)\eta(X) - da(X)]\eta(Y)\eta(Z)\\
\nonumber&+&(n-1)\eta(Z)\Omega(X,Y) - (n-3)\eta(Y)\Omega(X,Z)\\
\nonumber&+&(n-\frac{1}{2})\eta(Y)g(X,Z) - da(X)g(Y,Z).
\end{eqnarray}
By virtue of (\ref{eqn4.3}) and (\ref{eqn4.8}) we obtain from (\ref{eqn4.7}) that
\begin{eqnarray}
\label{eqn4.9}
(\nabla_{X}S)(Y,Z) &=& \eta(Z)S(Y,\phi X)-\eta(Y)[S(X,Z) - S(Z,\phi X)]\\
\nonumber&-&\eta(Y)\alpha(X,Z)-[(2n-1)\eta(X) - da(X)]\eta(Y)\eta(Z)\\
\nonumber&-&(n-1)\eta(Z)\Omega(X,Y) + (n-3)\eta(Y)\Omega(X,Z)\\
\nonumber&-&(n-\frac{1}{2})\eta(Y)g(X,Z) + da(X)g(Y,Z).
\end{eqnarray}
\indent Thus in a locally $\phi$-Ricci symmetric LP-Sasakian manifold with respect to semi-symmetric metric connection,
the relation (\ref{eqn4.9}) holds for any $X$, $Y$, $Z\in\chi(M)$.\\
\indent Next if the relation (\ref{eqn4.9}) holds in an LP-Sasakian manifold with respect to semi-symmetric metric connection
then for any horizontal vector fields $X$, $Y$, $Z$ with tr.$\alpha$ = constant, we obtain $(\nabla_{X}S)(Y,Z) = 0$
and hence the manifold is locally $\phi$-Ricci symmetric with respect to semi-symmetric metric connection.
Thus we can state the following:
\begin{thm}
An LP-Sasakian manifold $(M^n,g)$ is locally $\phi$-Ricci symmetric with respect to semi-symmetric metric connection
with tr.$\alpha$ = constant if and only if the relation \emph{(\ref{eqn4.9})} holds for any vector fields $X$, $Y$, $Z\in\chi(M)$.
\end{thm}
Putting $Y=\xi$ in (\ref{eqn4.9}) and using (\ref{eqn2.38}), we get
\begin{eqnarray}
\label{eqn4.10}
S(X,Z) &=& 2(n-2)\Omega(X,Z) - \alpha(X,Z)\\
\nonumber&-&(n-\frac{1}{2})g(X,Z) + (2n-1)\eta(X)\eta(Z)
\end{eqnarray}
for any vector fields $X$, $Z\in\chi(M)$.\\
This leads to the following:
\begin{thm}
In a locally $\phi$-Ricci symmetric LP-Sasakian manifold with respect to semi-symmetric metric connection,
the Ricci tensor is of the form \emph{(\ref{eqn4.10})}.
\end{thm}
\noindent{\bf Acknowledgements.} The authors wish to express their sincere
gratitude to the referee for the valuable suggestions to the improvement of the paper.

\vspace{0.1in}
\noindent S. K. Hui\\
Department of Mathematics\\
Sidho Kanho Birsha University\\
Purulia -- 723 104\\
West Bengal, India\\
E-mail: shyamal\_hui@yahoo.co.in\\

\noindent A. A. Shaikh\\
Department of Mathematics,\\
University of Burdwan,\\
Burdwan-713 104,\\
West Bengal, India\\
E-mail: aask2003@yahoo.co.in, aashaikh@math.buruniv.ac.in\\
\end{document}